\documentclass[11pt]{amsart} 
\usepackage{amscd}
\usepackage{amsfonts,amssymb,latexsym} % \Box, \twoheadrightarrow
\input xy
\xyoption{all}
\xyoption{arc}
\CompileMatrices

\newcommand{\SA}{{\mathcal{A}}}

\newcommand{\SH}{{\mathcal{H}}}
\newcommand{\SI}{{\mathcal{I}}}
\newcommand{\SJ}{{\mathcal{J}}}

\newcommand{\SO}{{\mathcal{O}}}
\newcommand{\SP}{{\mathcal{P}}}

\newcommand{\SU}{{\mathcal{U}}}
\newcommand{\SV}{{\mathcal{V}}}

\newcommand{\PP}{\mathbb{P}}
\newcommand{\ZZ}{\mathbb{Z}}
\newcommand{\CC}{\mathbb{C}}

\newcommand{\VV}{\mathbb{V}}

\newcommand{\FM}{\mathfrak{M}}

\newcommand{\isom}{\cong}

\newcommand{\Spec}{\operatorname{Spec}}

\newcommand{\codim}{\operatorname{codim}}

\newcommand{\Hom}{\operatorname{Hom}}

\newcommand{\Mor}{\operatorname{Mor}}
\newcommand{\Quot}{\operatorname{Quot}}

\newcommand{\Sch}{\operatorname{Sch}}
\newcommand{\Sets}{\operatorname{Sets}}
\newcommand{\Sym}{\operatorname{Sym}}
\newcommand{\id}{\operatorname{id}}

\newcommand{\surj}{\twoheadrightarrow}
\newcommand{\inj}{\hookrightarrow}
\newcommand{\too}{\longrightarrow}

\newcommand{\gr}{{\rm gr\,}}
\newcommand{\rk}{\operatorname{rk}}

\newcommand{\wt}{\widetilde}

\newcommand{\gl}{\operatorname{GL}}

\newcommand{\glv}{{\operatorname{GL}(V)}}

\newcommand{\slv}{\operatorname{SL}(V)}

\newcommand{\orth}{\operatorname{O}(r)}
\newcommand{\sor}{{\operatorname{SO}(r)}}
\newcommand{\sympl}{\operatorname{Sp}(r)}

\newcommand{\ev}{{\rm ev}}
\newcommand{\gitq}{/\!\!/}

\newcommand{\et}{{\rm et}}

\newcommand{\fg}{\mathfrak{g}}

\newcommand{\fgp}{{\mathfrak{\fg}'}}
\newcommand{\fz}{\mathfrak{z}}

\newcommand{\Aut}{{\rm Aut}}

\newcommand{\op}{\operatorname}

\newtheorem{proposition}{Proposition}[section]
\newtheorem{theorem}[proposition]{Theorem}
\newtheorem{definition}[proposition]{Definition}

\newtheorem{lemma}[proposition]{Lemma}

\numberwithin{equation}{section}

\title[Lectures on Principal bundles]
{Lectures on principal bundles over projective varieties}
\author[T. G\'omez]{Tom\'as L. G\'omez}
\address{
ICMAT, Consejo Superior de Investigaciones Cient{\'\i}ficas,
Serrano 113bis, 28006 Madrid, Spain; and 
Facultad de Matem\'aticas, 
Universidad Complutense de Madrid, 
28040 Madrid, Spain. }

\email{tomas.gomez@icmat.es}

\begin{document}

\begin{abstract}
Lectures given in the Mini-School
on Moduli Spaces at the Banach
center (Warsaw) 26-30 April 2005.
\end{abstract}
\maketitle

In these notes we will always work with schemes over the field of
complex numbers $\CC$.  Let $X$ be a scheme. A vector bundle of rank
$r$ on $X$ is a scheme with a surjective morphism $p:\VV\to X$ and an
equivalence class of linear atlases. A linear atlas is an open cover
$\{U_i\}$ of $X$ (in the Zariski topology) and isomorphisms
$\psi_i:p^{-1}(U_i) \to U_i\times \CC^r$, such that $p=p_X\circ
\psi_i$, and $\psi_j^{}\circ\psi_i^{-1}$ 
is linear on the fibers.  Two
atlases are equivalent if their union is an atlas.  These two
properties are usually expressed by saying that a vector bundle is
locally trivial (in the Zariski topology), and the fibers have a
linear structure.

An isomorphism of vector bundles on $X$ is an isomorphism $\varphi:\VV \to
\VV'$ 
of schemes which is compatible with the linear structure.  That
is, $p=p'\circ\varphi$ and the covering $\{U_i\}\bigcup \{U'_i\}$
together with the isomorphisms $\psi_i$, $\psi'_i\circ \varphi$ is a
linear structure on $\VV$ as before.

The set of isomorphism classes of vector bundles of rank $r$ on $X$ is
canonically bijective to the Cech cohomology set $\check
H^1(X,\underline{\gl_r})$.  
Indeed, since the transition functions
$\psi_j^{}\circ\psi_i^{-1}$ are linear on the fibers, they are given
by morphisms $\alpha_{ij}:U_i\cap U_j \to \gl_r$ which satisfy the
cocycle condition.

Given a vector bundle $\VV\to X$ we define the locally free sheaf $E$
of its sections, which to each open subset $U\subset X$, assigns
$E(U)=\Gamma(U,p^{-1}(U))$.  This provides an equivalence of
categories between the categories of vector bundles and that of locally
free sheaves (\cite[Ex. II.5.18]{Ha}).  Therefore, if no confusion
seems likely to arise, we will use the words ``vector bundle'' and
``locally free sheaf'' interchangeably.  Note that a vector bundle of
rank 1 is the same thing as a line bundle.  We will be interested in
constructing moduli spaces of vector bundles, which can then be
considered as generalizations of the Jacobian.  Sometimes it will be
necessary to consider also torsion free sheaves.  For instance, in
order to compactify the moduli space of vector bundles.

Let $X$ be a smooth projective variety of dimension $n$ with an ample
line bundle $\SO_X(1)$ corresponding to a divisor $H$.  Let $E$ be a
torsion free sheaf on $X$.  Its Chern classes are denoted $c_i(E)\in
H^{2i}(X;\CC)$.  
We define the \textit{degree} of $E$
$$
\deg E = c_1(E) H^{n-1}
$$ 
and its Hilbert polynomial
$$
P_E(m)=\chi(E(m)),
$$
where $E(m)=E\otimes \SO_X(m)$ and $\SO_X(m)=\SO_X(1)^{\otimes m}$.

If $E$ is locally free, we define the determinant line bundle as
$\det E=\bigwedge^r E$. If $E$ is torsion free,
since $X$ is smooth, we can still define its determinant as follows.
The maximal open subset $U\subset X$ where
$E$ is locally free is \textit{big} (with this we will mean that its
complement has codimension at least two), because it is torsion free.
Therefore, there is a line bundle $\det E|_U$ on $U$, and since
$U$ is big and $X$ is smooth, this extends to a unique line bundle
on $X$, which we call the determinant $\det E$ of $E$.
It can be proved that $\deg E=\deg \det E$.

We will use the following notation. Whenever ``(semi)stable''
and  ``$(\leq)$'' appears in a sentence, two sentences should
be read. One with ``semistable'' and ``$\leq$'', and another
with ``stable'' and ``$<$''. Given two polynomials $p$ and $q$,
we write $p<q$ if $p(m)<q(m)$ when $m\gg 0$.

A torsion free sheaf $E$ is \textit{(semi)stable} if for all proper
subsheaves $F\subset E$, 
$$
\frac{P_F}{\rk F} (\leq) \frac{P_E}{\rk E} \; .
$$
A sheaf is called \textit{unstable} if it is not semistable.
Sometimes this is referred to as Gieseker (or Maruyama) stability.

A torsion free sheaf $E$ is \textit{slope-(semi)stable} if for all proper
subsheaves $F\subset E$ with $\rk F<\rk E$, 
$$
\frac{\deg{F}}{\rk F} (\leq) \frac{\deg{E}}{\rk E} \; .
$$
The number $\deg E/\rk E$ is called the \textit{slope} of $E$.
A sheaf is called \textit{slope-unstable} if it is not slope-semistable.
Sometimes this is referred to as Mumford (or Takemoto) stability.

Using Riemann-Roch theorem, we find
$$
P_E(m)= \rk E \frac{m^n}{n!} + (\deg E - \rk E \frac{\deg K}{2})
\frac{m^{n-1}}{(n-1)!} + \cdots
$$
where $K$ is the canonical divisor.
{}From this it follows that
$$
\text{slope-stable $\Longrightarrow$ stable
$\Longrightarrow$ semistable $\Longrightarrow$ 
slope-semistable}
$$
Note that, if $n=1$, Gieseker and Mumford (semi)stability coincide,
because the Hilbert polynomial has degree 1.

Let $E$ be a torsion free sheaf on $X$. There is a unique filtration,
called the Harder-Narasimhan filtration,
$$
0 = E_0 \subsetneq E_1 \subsetneq E_2 \subsetneq \cdots \subsetneq E_l=E
$$
such that $E^i=E_i/E_{i-1}$ is semistable and
$$
\frac{P_{E^i}}{\rk E^i} > \frac{P_{E^{i+1}}}{\rk E^{i+1}} 
$$
for all $i$. In particular, any torsion free sheaf can be described
as successive extensions of semistable sheaves.

There is also a Harder-Narasimhan filtration for slope stability: this is
the unique filtration such that $E^i=E_i/E_{i-1}$ is slope-semistable and
$$
\frac{\deg {E^i}}{\rk E^i} > \frac{\deg {E^{i+1}}}{\rk E^{i+1}} 
$$
for all $i$. We denote
$$
\mu_{\max}(E)=\mu(E^1)\; ,\quad  \mu_{\min}(E)=\mu(E^l).
$$
Of course, in general these two filtrations will be different.

Now let $E$ be a semistable sheaf. There is a filtration, called
the Jordan-H\"older filtration,
$$
0 = E_0 \subsetneq E_1 \subsetneq E_2 \subsetneq \cdots \subsetneq E_l=E
$$
such that $E^i=E_i/E_{i-1}$ is stable and
$$
\frac{P_{E^i}}{\rk E^i} = \frac{P_{E^{i+1}}}{\rk E^{i+1}} 
$$
for all $i$. This filtration is not unique, but the associated
graded sheaf
$$
\gr_{\rm JH}(E)=\bigoplus_{i=1}^l E^i
$$
is unique up to isomorphism. It is easy to check that $\gr_{\rm JH}(E)$
is semistable.
Two semistable torsion free sheaves are
called \textit{$S$-equivalent} if $\gr_{\rm JH}(E)$ and $\gr_{\rm JH}(E')$
are isomorphic. 

There is also a Jordan-H\"older filtration for slope stability, just 
replacing Hilbert polynomials with degrees.

A \textit{family of coherent sheaves} parameterized by a scheme
$T$ (also called $T$-family) 
is a coherent sheaf $E_T$ on $X\times T$, flat over $T$. 
For each closed
point $t\in T$, we get a sheaf $E_t:=f^* E_T$ on $X\times t \cong X$,
where $f:X\times t \to X\times T$ is the natural inclusion. 
We say that $E_T$ is a family of torsion-free sheaves
if $E_t$ is torsion-free sheaf  
for all closed points $t\in T$. 
We have analogous definitions for any open condition, and
hence we can talk of families of (semi)stable sheaves, 
of families of sheaves with fixed Chern classes $c_i(E)$,
etc...
Two families are \textit{isomorphic} if $E^{}_T$ and $E'_T$ are
isomorphic as sheaves.

To define the notion of moduli space, we will first look at
the Jacobian $J$ of a projective scheme $X$. 
There is a bijection between isomorphism classes of line bundles
$L$ with $0=c_1(L)\in H^2(X;\CC)$ and closed points of $J$. 

Furthermore, if we are given a family of line bundles $L_T$, with
vanishing first Chern class, parameterized by a scheme $T$, 
we obtain a morphism $f:T\to J$ such that for all $t\in T$, 
the point $f(t)\in J$
is the point corresponding to the isomorphism class of $L_t$.
And, conversely, if we are given a morphism $f:T\to J$, we obtain
a family of line bundles parameterized by $T$ by pulling-back
a Poincare line bundle: $L_T =(\id_X\times f)^* \SP$.

Note that both constructions are not quite inverse to each other.
On the one hand, if $M$ is a line bundle on $T$, the families
$L_T$ and $L_T \otimes p_T^* M$ give the same morphism from 
$T$ to the Jacobian, and on the other hand, there is no unique
Poincare line bundle: given a linen bundle $M$ on $J$, 
$\SP\otimes p_J^* M$ is also a Poincare line bundle, 
and the family induced by $f$ will change to $L_T \otimes 
p^*_T f^*_{}M$. 
This is why we declare two families of line bundles
\textit{equivalent} if they differ by the pullback of a line bundle
on the parameter space $T$. 

Using this equivalence, both constructions become inverse of each
other. That is, there is a bijection between equivalence
classes of families and morphisms to the Jacobian.

One could ask: is there a ``better'' version of the Jacobian?, i.e.
some object $\SJ$ which provides a bijection between morphisms to it
and families of line bundles up to isomorphism (not up to equivalence).
The answer is yes, but this object $\SJ$ is not a scheme!. It is
an algebraic stack (in the sense of Artin): the Jacobian stack. 
A stack is a generalization of the notion of scheme, but we will not
consider it here.

We would like to have a scheme with the same properties as the Jacobian, 
but for torsion-free sheaves instead of line bundles.
To be able to do this, we have to consider only the semistable ones.
Then there will be a moduli scheme $\FM(r,c_i)$ such that a family of
semistable torsion-free sheaves parameterized by $T$, 
with rank $r$ and Chern classes
$c_i$, will induce a morphism from the
parameter space $T$ to $\FM(r,c_i)$. In particular, to each semistable torsion
free sheaf we associate a closed point. 
If two stable torsion free sheaves on $X$ are not isomorphic,
they will correspond to different points of $\FM(r,c_i)$, but
it can happen that two strictly semistable torsion free sheaves on $X$ 
which are not isomorphic correspond to
the same point of $\FM(r,c_i)$. 
In fact, $E$ and $E'$ correspond to the same point if and
only if they are S-equivalent.

Another difference with the properties of the Jacobian is that
in general there will be no
``universal torsion free sheaf'' on $X\times \FM(r,c_i)$, 
i.e. there will be no analogue
of the Poincare bundle. In other words, given a morphism 
$f:T\to \FM(r,c_i)$, there might be no family parameterized by $T$ which
induces $f$. If there is a universal torsion-free sheaf, we say
that $\FM(r,c_i)$ is a \textit{fine moduli space}, and if it does not exist,
we say that it is a \textit{coarse moduli space}.

To explain this more precisely, it is useful to use the language of
representable functors. Given a scheme $M$ over $\CC$, we define
a (contravariant) functor $\underline{M}:=\Mor(-,M)$ from the category of
$\CC$-schemes $(\Sch/\CC)$ to the category of sets $(\Sets)$ by sending
an $\CC$-scheme $B$ to the set of morphisms $\Mor(B,M)$. On morphisms
it is defined with composition, i.e., to a morphism $f:B\to B'$ we
associate the map $\Mor(B',M)\to \Mor(B,M)$ which sends $\varphi'$ 
to $\varphi\circ f$. 

\begin{definition}[Represents]
A functor $F:(\Sch/\CC) \to (\Sets)$ is represented by a scheme $M$
if there is an isomorphism of functors $F\cong \underline{M}$.
\end{definition}

Of course, not all functors
from $(\Sch/\CC)$ to $(\Sets)$
are representable, but if a functor $F$ is, then the scheme $M$ is 
unique up to canonical isomorphism. Given a morphism $f:M\to M'$, we
obtain an natural transformation $\underline{M}\to \underline{M'}$,
and, by Yoneda's lemma, every natural transformation between 
representable functors is induced by a morphism of schemes. In other
words, the category of schemes is a full subcategory of the category
of functors $(\Sch/\CC)'$, whose objects are contravariant functors from
$(\Sch/\CC)$ to $(\Sets)$ and whose morphisms are natural transformation.
Therefore, we will denote by the same letter a morphism of schemes and
the associated natural transformation.

For instance, let $F_J:(\Sch/\CC)\to (\Sets)$ be the functor which 
sends a scheme $T$ to the set of equivalence classes of $T$-families of line
bundles on $X$, with $c_1=0$.
This functor is represented by the
Jacobian, i.e., there is an isomorphism of functors $F_J\cong
\underline{J}$. 
This is the translation, to the language of
representable functors, of the fact that there is a 
natural bijection between
the set of equivalence classes of these families and the set of
morphisms from $T$ to $J$.

\begin{definition}[Corepresents]
A functor $F:(\Sch/\CC) \to (\Sets)$ is corepresented by a scheme 
$M$ if there is a natural transformation of functors $\phi:
F \to \underline{M}$ such that
given another scheme $M'$ and natural transformation $\phi':F
\to \underline{M'}$, there is a unique morphism $\eta:
M\to M'$ with $\phi'= \eta \circ \phi$.
$$
\xymatrix{
{F} \ar[d]_{\phi} \ar[rd]^{\phi'} \\
 {\underline{M}} \ar@{-->}[r]_{\exists !\, \eta} &  {\underline{M'}}\\ 
}
$$
\end{definition}

If $M$ corepresents $F$, then $M$ is unique up to canonical isomorphism.
To explain why this is called ``corepresentation'', 
let $(\Sch/\CC)'$ be the above defined
functor category. Then it can be seen that
$M$ represents $F$ if and only 
if there is a natural bijection $\Mor(Y,M)=
\Mor_{(\Sch/\CC)'}(\underline{Y},F)$ for all schemes $Y$.
On the other hand,
$M$ corepresents $F$ if and only if there is a natural 
bijection $\Mor(M,Y)=
\Mor_{(\Sch/\CC)'}(F,\underline{Y})$ for all schemes $Y$.
If $M$ represents $F$, then it corepresents it, but the converse is
not true.

Let $X$ be a fixed $\CC$-scheme. 
Define a functor $F^{ss}_{r,c_i}$ from the category of schemes over $\CC$ to
the category of sets, sending a scheme $T$ to the set 
$F^{ss}_{r,c_i}(T)$
of isomorphism classes of families of torsion-free sheaves 
on $X$
parameterized by $T$, with rank $r$ and Chern classes $c_i$.
On morphisms it is defined by pullback, i.e., to a morphism
$f:T\to T'$ we associate the map $F(T')\to F(T)$ 
which sends the family $E_T'$ to $(\id_X \times f )^*E'_T$.
Analogously, we define the functor $F^{s}_{r,c_i}$ of families of
stable torsion free sheaves. 

It can be shown that, for any polarized smooth projective variety $X$,
there is a scheme $\FM(r,c_i)$ corepresenting the above defined
functor $F^{ss}_{r,c_i}$ (\cite{Gi,Ma,Sesh,Si}). In section \ref{secsheaves}
we will sketch a proof of this result.

Note that the transformation of functors
$\phi$ gives, for any $T$-family of semistable torsion free sheaves,
a morphism $f:T\to \FM(r,c_i)$. 
As we mentioned before, there is a canonical bijection between
closed points of $\FM(r,c_i)$ and S-equivalence classes of
semistable torsion free sheaves.

Let $\hat F^{ss}_{r,c_i}$ be the functor of equivalence classes of
families of semistable sheaves, where, as before,
we declare two families equivalent
if they differ by the pullback of a line bundle on the parameter
space. There are some cases in which this functor is representable
(for instance, if the rank $r$ and degree $c_1$ are coprime). In these
cases, there is a universal family parameterized by the moduli space,
and this universal family is unique up to the pullback of a line
bundle on the moduli space. 

\begin{definition}[Moduli space]
We say that $M$ is a moduli space for a set of objects, if 
it corepresents the functor of families of those objects.
\end{definition}

\begin{definition}[Coarse moduli]
A scheme $M$ is called a coarse moduli scheme for $F$ 
if it corepresents $F$
and furthermore
the map 
$$
\phi(\Spec \CC):F(\Spec \CC) \to \Hom(\Spec \CC, M)
$$ 
is bijective.
\end{definition}

Note that if a functor $F$ is corepresented by a scheme $M$, then
it is a coarse moduli scheme for the functor $\wt F$ of S-equivalence
classes of $F$, i.e., the functor defined as 
$$
\wt F(T)=
\left\{
\begin{array}{ll}
F(T) & \text{, if $T\neq \Spec \CC$}\\
\text{S-equivalence classes of objects of $F(\Spec \CC)$}
& \text{, if $T= \Spec \CC$}\\
\end{array}
\right .
$$

\section{Moduli space of torsion free sheaves}
\label{secsheaves}

In this section we will sketch the proof of the existence of 
the moduli space of semistable torsion free sheaves.
We will start by giving a brief 
idea of the construction.
It can be shown
that there is a scheme $Y$ classifying 
\textit{semistable based sheaves}, that is,
pairs $(f,E)$, where
$E$ is a semistable sheaf and $f:V \to H^0(E(m))$ is an isomorphism
between a fixed vector space $V$ and $H^0(E(m))$. The group $\slv$ 
acts on $Y$ by ``base change'': an element $g\in \slv$ sends the pair
$(f,E)$ to $(f\circ g,E)$. Two pairs $(f,E)$ and $(f', E')$ are
in the same orbit if and only if $E$ is isomorphic to $E'$, therefore, 
the quotient of $Y$ by the action of $\slv$ 
will be a moduli space of semistable sheaves.
But, does this quotient exist in the category of schemes?, i.e., is there
a scheme whose points are in bijection with the $\slv$-orbits in $Y$?.
In general the answer is no, but Geometric Invariant Theory (GIT) gives us
something which is quite close to this, and is called the GIT
quotient, 
and this will be the moduli space.

Note that we are using the group $\slv$, and not $\glv$. This is because
if two isomorphisms $f$ and $f'$ only differ by multiplication with a 
scalar, then they correspond to the same point in $Y$. In other words,
$Y$ classifies pairs $(f,E)$ up to scalar.

Let $G$ be an algebraic group. Recall that a 
right action on a scheme $R$ is a 
morphism $\sigma:R\times G\to R$, which we will usually denote
$\sigma(z,g)= z\cdot g$,
such that $z\cdot (gh)=(z\cdot g)\cdot h$ and $z\cdot e=z$, where
$e$ is the identity element of $G$. A left action is analogously
defined, with the associative condition $(hg)\cdot z= h\cdot (g\cdot z)$.

The orbit of a point $z\in R$
is the image $z \cdot G$. A morphism $p:R \to M$ between two schemes
endowed with $G$-actions is called \textit{$G$-equivariant} if it commutes
with the actions, that is $f(z)\cdot g= f(z\cdot g)$. If the action
on $M$ is trivial (i.e. $y\cdot g=y$ for all $g\in G$ and $y\in M$), 
then we also 
say that $f$ is \textit{$G$-invariant}.

If $G$ acts on a projective variety $R$, a \textit{linearization} 
of the action on a line bundle $\SO_R(1)$ 
consists of giving, for each $g\in G$,
an isomorphism of line bundles 
$\wt g:\SO_R(1)\to \varphi_g^* \SO_R(1)$,
($\varphi_g=\sigma(\cdot,g)$)
which also satisfies the
previous associative property. Giving a linearization is thus the same
thing as giving an
action on the total space $\VV$ of the line bundle, which is linear 
along the fibers, and such that the projection $\VV\to R$ is
equivariant. If $\SO_R(1)$ is very ample, then a linearization
is the same thing as a representation of $G$ on the vector space
$H^0(\SO_R(1))$ 
such that the natural embedding $R\to \PP(H^0(\SO_R(1))^\vee)$
is equivariant.

\begin{definition}[Categorical quotient]
Let $R$ be a scheme endowed with a $G$-action. 
A categorical
quotient is a scheme $M$ with a $G$-invariant morphism $p:R\to M$ 
such that for every other scheme $M'$, and $G$-invariant morphism
$p'$, there is a unique morphism $\varphi$ with $p'=\varphi\circ p$
$$
\xymatrix{
{R} \ar_{p}[d] \ar^{p'}[rd]\\
%{M} \ar_{\exists !\, \varphi}[r]& {M'}
{M} \ar@{-->}_{\exists !\, \varphi}[r]& {M'}
}
$$
\end{definition}

\begin{definition}[Good quotient]
Let $R$ be a scheme endowed with a $G$-action. 
A good quotient 
is a scheme $M$ with a $G$-invariant morphism $p:R\to M$ 
such that
\begin{enumerate}
\item $p$ is surjective and affine

\item $p_*(\SO^G_R)=\SO_M^{}$, where $\SO^G_R$ is the sheaf of
$G$-invariant functions on $R$.

\item If $Z$ is a closed $G$-invariant subset of $R$, then $p(Z)$ is
closed in $M$. Furthermore, if $Z_1$ and $Z_2$ are two 
closed $G$-invariant subsets of $R$ with $Z_1\cap Z_2=\emptyset$, 
then $f(Z_1)\cap f(Z_2)=\emptyset$.
\end{enumerate}
\end{definition}

\begin{definition}[Geometric quotient]
A geometric quotient $p:R\to M$ is a good quotient 
such that $p(x_1)=p(x_2)$ if and only if
the orbit of $x_1$ is equal to the orbit of $x_2$.
\end{definition}

Clearly, a geometric quotient is a good quotient, and 
a good quotient is a categorical quotient.

Geometric Invariant Theory (GIT) is a technique to construct good
quotients (cf. \cite{Mu1}). 
Assume $R$ is projective, and the action of $G$ on $R$ has
a linearization on an ample line bundle $\SO_R(1)$.
A closed point $z\in R$ is called \textit{GIT-semistable} 
if, for some
$m>0$, there is a $G$-invariant section $s$ of $\SO_R(m)$
such that $s(z)\neq 0$. 
If, moreover, 
the orbit of $z$ is closed in the open set of all
GIT-semistable points, it is called \textit{GIT-polystable},
and, if furthermore, this closed orbit has 
the same dimension as $G$ ( i.e., if $z$ has finite stabilizer),
then $z$ is called a \textit{GIT-stable} point.
We say that a closed point of $R$ is \textit{GIT-unstable}
if it is not GIT-semistable.

We will use the following characterization in 
\cite{Mu1} of GIT-(semi)stability. Let 
$\lambda:\CC^*\to G$ be a one-parameter subgroup
(by this we mean a nontrivial group homomorphism,
even if $\lambda$ is not injective),
and let $z\in R$. Then $\lim_{t\to 0} z\cdot \lambda(t)= z_0$
exists, and $z_0$ is fixed by $\lambda$. Let
$t\mapsto t^a$ be the character by which $\lambda$
acts on the fiber of $\SO_R(1)$. Defining 
$\mu(z,\lambda)=a$, Mumford proves that 
$z$ is GIT-(semi)stable if and
only if, for all one-parameter subgroups, 
it is $\mu(z,\lambda)(\leq)0$.

\begin{proposition}
Let $R^{ss}$ (respectively, $R^{s}$) be the subset of
GIT-semistable points (respectively, GIT-stable).
Both $R^{ss}$ and $R^{s}$ are open subsets. 
There is a good quotient $R^{ss}\to R\gitq G$,
the image $R^{s}\gitq G$ of $R^{s}$ is open, 
$R\gitq G$ is projective,
and the restriction $R^{s}\to R^{s}\gitq G$ 
is a geometric quotient.
\end{proposition}

There is one important case in which a scheme is only 
quasi-projective but GIT can be applied to get a projective
quotient:
Assume that $R'$ is a $G$-acted scheme with a linearization 
on a line bundle $\SO_{R'}(1)$, which is the restriction of 
a linearization on an ample line bundle $\SO_R(1)$ on a projective
variety $R$, and $R'=R^{ss}$, the open subset of GIT-semistable points
of $R$. Then we define $R'\gitq G=R\gitq G$.

Now we are going to describe Grothendieck's Quot-scheme.
This scheme parameterizes quotients of a fixed coherent sheaf
$\SV$ on $X$. That is, pairs $(q, E)$, where $q:\SV\surj E$
is a surjective homomorphism and $E$ is a coherent sheaf on $X$.
An isomorphism of quotients is an isomorphism $\alpha:E\to E'$ such
that the following diagram is commutative
\[
\xymatrix{
{\SV} \ar@{=}[d] \ar@{->>}^{q}[r] & {E}\ar[d]_{\isom}^{\alpha}\\
{\SV} \ar@{->>}^{q'}[r] & {E'}\\
}
\]
A family of quotients parameterized by $T$ 
is a pair $(q:p_X^*\SV \surj
E_T,E_T)$ where $E_T$ is a coherent sheaf on $X\times T$, flat 
over $T$. An isomorphism of families is an isomorphism $\alpha:E^{}_T
\to E'_T$ 
such that $\alpha\circ q=q'$.
Recall that $X$ is a projective scheme endowed with an ample line
bundle $\SO_X(1)$. Therefore, if $E_T$ is flat
over $T$ then the Hilbert polynomial $P_{E_t}$ is
locally constant as a function of $t\in T$. If $T$ is reduced, the
converse is also true.

Fix a polynomial $P$ and a coherent sheaf $\SV$ on $X$. Consider the
contravariant functor which sends a scheme $T$ to the set of
isomorphism classes of $T$-families of sheaves with Hilbert polynomial
$P$ (and it is defined as pullback on morphisms). Grothendieck proved
that there is a projective scheme $\Quot_X(\SV,P)$, called the Quot
scheme, which represents this functor. In particular, there is
a universal quotient, i.e.,
a tautological family of quotients parameterized by $\Quot_X(\SV,P)$.
We will be interested in the case $\SV=V\otimes_\CC \SO_X(-m)$, where
$V$ is a vector space and $m$ is sufficiently large. 

Given a coherent sheaf $E$, 
there is an integer $m(E)$ 
such that, if $m\geq m(E)$,
then $E(m)$ is generated by global sections,
$h^0(E(m))=P_E(m)$, and 
$h^i(E(m))=0$ for $i>0$
(\cite[Def 1.7.1]{H-L}).
Assume that $m\geq m(E)$ and $\dim V=P_E(m)$.
An isomorphism $f:V \to H^0(E(m))$ induces a quotient
$$
q:V\otimes_\CC \SO_X(-m) \stackrel{\cong}\too 
H^0(E(m))\otimes_\CC\SO_X(-m)
\too E
$$
as above, and this is how a scheme parameterizing based sheaves $(f,E)$
appears as a subscheme of Grothendieck's Quot scheme.

Note that if we have a set $\SA$ of isomorphism classes of sheaves, 
there might not be an integer $m$ large enough for all sheaves.
A set $\SA$ of isomorphism classes of sheaves on $X$ is called
\textit{bounded} if there is a family
$E_S$ of torsion free sheaves parameterized by a scheme
$S$ of finite type, such that for all $E\in \SA$, there is at least one
point $s\in S$ such that the corresponding sheaf $E_s$ 
is isomorphic to $E$.
If a set $\SA$ is bounded, then we can find an integer $m$
such that $m\geq m(E)$ for all $E\in \SA$, thanks to the fact that $S$
is of finite type.

Maruyama proved that 
the set $\SA$ of semistable sheaves with fixed Hilbert polynomial 
is bounded, and it follows that
there is an integer $m_0$, depending only on the polynomial $P$ 
and $(X,\SO_X(1))$, such that $m_0\geq m(E)$ for all semistable
sheaves $E$.
This technical result is crucial in order to construct the moduli
space. In fact, if $\dim X>1$, he was able to prove it
only if the base field has characteristic 0, and therefore he
could only prove the existence of the moduli space in this case.
Recently Langer was able to prove boundedness for characteristic
$p>0$, and therefore he was able to construct the corresponding
moduli space \cite{La}. 

Fix a Hilbert polynomial $P$, and let $m\geq m_0$.
Let $Y\subset \Quot_X(V\otimes_\CC \SO_X(-m),P)$ be the 
open subset of quotients such that $E$ is torsion free
and $q$ induces an isomorphism $V\cong H^0(E(m))$.
Let  $\overline{Y}$ be the closure of the open set 
$Y$ in $\Quot_X(V\otimes_\CC \SO_X(-m),P)$.
Note that there is a natural action of $\slv$ on 
$\Quot_X(V\otimes_\CC \SO_X(-m))$, 
which sends the quotient $q:V\otimes_\CC\SO_X(-m) \to E$
to the composition $q\circ (g\times \id)$. It leaves $Y$ and $\overline{Y}$ 
invariant, and coincides with the previously defined action for 
based sheaves $(f,E)$. 

To apply GIT, we also need an ample line bundle 
on $\overline{Y}$ and a linearization of the $\slv$-action on it.
This is done by giving an embedding of $\overline{Y}$ in $\PP(V_1)$, 
where $V_1$ will be a vector space with a
representation of $\slv$.

There are different ways of doing this, corresponding to different
representations $V_1$.
One of them corresponds to Grothendieck's embedding of 
the Quot scheme. This is the method used by Simpson
\cite{Si}. Let $q:V\otimes_\CC \SO_X(-m)\surj E$ be a quotient.
Let $l>m$ be an integer and $W=H^0(\SO_X(l-m))$.
The  quotient $q$ 
induces homomorphisms
$$
\begin{array}{rccl}
q\,\,: & V\otimes_\CC \SO_X(l-m)&\surj& E(l) \\
q'\,:& V\otimes W&\to& H^0(E(l)) \\
q'':& \bigwedge{}^{P(l)}(V\otimes W) &\to& 
\bigwedge{}^{P(l)} H^0(E(l))\;\cong \; \CC
\end{array}
$$
If $l$ is large enough, these homomorphisms are surjective, and
give Grothendieck's embedding of the Quot scheme.
$$
\Quot_X(V\otimes_\CC \SO_X(-m),P) \;\too\; 
\PP\Big(\bigwedge{}^{P(l)}(V^\vee\otimes W^\vee)\Big),
$$
The natural representation of $\slv$ in
$\bigwedge{}^{P(l)}(V^\vee\otimes W^\vee)$  
gives a linearization of the $\slv$ action on the
very ample line bundle $\SO_{\overline{Y}}(1)$ induced by this 
embedding on $\overline{Y}$.

A theorem of Simpson says that a point $(q,E)\in \overline{Y}$ 
is GIT-(semi)stable
if and only if the sheaf $E$ is (semi)stable and the induced linear map
$f:V \to H^0(E(m))$ is an isomorphism. In other words, $Y=\overline{Y}^{ss}$.
Therefore, 
the GIT quotient $\overline{Y}\gitq \slv$ is the moduli space $\FM(P)$ 
of semistable
sheaves with Hilbert polynomial $P$. The Chern classes $c_i\in
H^{2i}(X,\CC)$
in a family of sheaves are locally constant, therefore the moduli
space $\FM(r,c_i)$ of semistable sheaves with fixed rank and Chern
classes is a union of connected components of the scheme $\FM(P)$.

Another choice of representation $V_1$ (and therefore, of line
bundle on $Y$ and linearization of the action) is the one
used by Gieseker and Maruyama. It is explained in the lectures of 
Schmitt.

\section{Moduli space of tensors}
\label{sectensors}

A \textit{tensor} of type $a$ is a pair $(E,\varphi)$ where
$E$ is a torsion free sheaf and 
$$
\varphi: E  \overbrace{\otimes \cdots \otimes}^a E \too \SO_X
$$
is a homomorphism.
An isomorphism between the tensors $(E,\varphi)$ and $(E',\varphi')$ is
a pair $(f,\alpha)$ where $f$ is an isomorphism between $E$ and $E'$,
$\alpha\in \CC^*$, and the following diagram commutes
$$
\xymatrix{
{E^{\otimes a}} \ar[r]^{\varphi} \ar[d]_{f^{\otimes a}} & 
{\SO_X} \ar[d]^{\alpha} \\
{E'{}^{\otimes a}} \ar[r]^{\varphi'} & {\SO_X} \\
}
$$
The definition of families of tensors and their isomorphisms are left
to the reader (\cite{G-S1, GLSS}).

To define the notion of stability for tensors, it is not enough to
look at subsheaves. We have to consider filtrations $E_\bullet\subset E$.
By this we always
understand a $\ZZ$-indexed filtration
$$
\ldots \subset E_{i-1}\subset E_{i}
\subset E_{i+1}\subset \ldots
$$
starting with $0$ and ending with $E$ (i.e., $E_k=0$ and $E_l=E$ for
some $k$ and $l$). 
We say that the filtration is 
\textit{saturated} if $E^i=E_i/E_{i-1}$ is torsion free for all $i$.
If we delete, from $0$ onward, all the non-strict 
inclusions, we obtain a filtration 
$$
0 \subsetneq E_{\lambda_1} \subsetneq E_{\lambda_2} 
\subsetneq \;\cdots\; \subsetneq E_{\lambda_t} \subsetneq
E_{\lambda_{t+1}}=E \qquad \lambda_1<\cdots<\lambda_{t+1}
$$
Reciprocally,
from a filtration $E_{\lambda_\bullet}$ we recover the
$\ZZ$-indexed filtration
$E_\bullet$ by 
defining $E_m=E_{\lambda_{i(m)}}$, 
where $i(m)$ is the maximum 
index with $\lambda_{i(m)}\leq m$.

Let $\SI_a=\{1,\ldots,t+1\}^{\times a}$ be the set of all
multi-indexes $I=(i_1,\ldots,i_a)$ of cardinality $a$. 
Define
\begin{equation}
\label{muE}	
\mu_{\rm tens}(\varphi,E_{\lambda_\bullet})=\min_{I\in \SI_a} \big\{
\lambda_{i_1}+\dots+\lambda_{i_a}: \,
\phi|^{}_{E_{\lambda_{i_1}}\otimes\cdots \otimes E_{\lambda_{i_a}}}\neq 0
  \big\} \; ,
\end{equation}
or, in terms of the $\ZZ$-indexed filtration,
\begin{equation}
\mu_{\rm tens}(\varphi,E_{\bullet})=\min_{I\in \SI_a} \big\{
{i_1}+\dots+{i_a}: \,
\phi|^{}_{E_{i_1}\otimes\cdots \otimes E_{i_a}}\neq 0
  \big\}
\end{equation}

\begin{definition}[Balanced filtration]
\label{defbalancedfiltration}
A saturated filtration $E_{\bullet}\subset E$ 
of a torsion free sheaf $E$ is called a balanced filtration 
if $\sum i \rk{E^i}=0$ for $E^i=E_i/E_{i-1}$.
In terms of $E_{\lambda_\bullet}$,
this is $\sum_{i=1}^{t+1} \lambda_i \rk (E^{\lambda_i})=0$
for $E^{\lambda_i}=E_{\lambda_i}/E_{\lambda_{i-1}}$.
\end{definition}

\begin{definition}[Stability of tensors]
Let $\delta$ be a polynomial of degree at most $n -1$ 
(recall $n=\dim X$) with positive
leading coefficient. 
We say that a tensor $(E,\varphi)$ is $\delta$-(semi)stable if 
$\varphi$ is not identically zero and
for all balanced
filtrations $E_{\lambda_\bullet}$ of $E$, it is 
\begin{equation}
\label{stabtensor}
\Big(\sum_{i=1}^t (\lambda_{i+1}-\lambda_i)
\big( r P_{E_{\lambda_i}} -r_{\lambda_i} P \big)\Big) +
\mu_{\rm tens}(\phi,E_{\lambda_\bullet}) \, \delta \;(\leq)\; 0
\end{equation}
\end{definition}

We will always denote $r=\rk E$ and $r_i=\rk E_i$.
The notion of stability for tensors looks complicated, but one finds that,
in the applications,
when the tensor has some geometric meaning, it can be simplified.
We will see some examples.

A \textit{framed bundle} is a tensor
of the form $(E,\varphi:E \to \SO_X)$. If $E$ is a vector bundle, then 
taking the dual we have a section of $E^\vee$, so this is equivalent
to the pairs $(E,\varphi: \SO_X\to E)$ considered by Bradlow,
Garc{\'\i}a-Prada and others. In this case, it is enough to look
at filtrations with one step, i.e. subsheaves $E'\subsetneq E$.

An \textit{orthogonal sheaf} 
is a tensor of the form $(E, \varphi:E \otimes E \to \SO_X)$, where 
$E$ is torsion free and $\varphi$
is symmetric and non-degenerate (in the sense that the induced
homomorphism $\det E\to \det E^\vee$ is an isomorphism).
A \textit{symplectic sheaf} is analogously
defined, requiring the tensor $\varphi$ to be skew-symmetric
instead of symmetric.

Given a subsheaf $E'\subset E$, its orthogonal
$E'{}^\perp$ is defined as the kernel of the composition
$$
E\stackrel{\wt\varphi}\too E^\vee \too E'{}^{\perp}\; ,
$$
where $\wt\varphi$ is induced by $\varphi$.

\begin{definition}
An orthogonal (or symplectic) sheaf is (semi)stable if for all
orthogonal filtrations, that is, filtrations with
$$
E_i^\perp = E_{-i-1}  
$$
for all $i$, the following holds
$$
\sum (r P_{E_i} - r_i P_E ) (\leq) 0 \, .
$$
\end{definition}

It it shown in \cite{G-S1} that an orthogonal (or symplectic) sheaf is
(semi)stable if and only if it is $\delta$-(semi)stable as a tensor,
when $\delta$ has degree $n-1$. 

A $T$-family of orthogonal sheaves is a $T$-family of tensors
$(E_T,\varphi_T:E_T\otimes E_T\too \SO_{X\times T})$ such that
$\varphi_T$ is symmetric and non-degenerate. Note that, since
being symmetric is a closed condition, it is not enough to check
that $\varphi_t$ is symmetric for every point $t\in T$. On the
other hand, being non-degenerate is an open condition, so it is
enough to check it for $\varphi_t$, for all points $t\in T$.

A \textit{Lie algebra sheaf} is a pair $(E,\varphi)$ where $E$ is a torsion
free sheaf and 
$$
\varphi: E \otimes E \too E^{\vee\vee}
$$
is a homomorphism such that for each point $x\in X$, where $E$ is
locally free, the induced homomorphism  on the fiber 
$\varphi(x):E(x)\otimes E(x) \to E(x)$ is a Lie algebra structure.
An isomorphism to another Lie algebra sheaf $(E',\varphi')$ 
is an isomorphism of
sheaves $f:E\to E'$ with $\varphi'\circ (f\otimes f)=f\circ \varphi$.

At first sight, this does not seem to be included in the formalism of
tensors, but, using the canonical
isomorphism 
\begin{equation}
\label{wedge}
(\bigwedge^{r-1} E)^\vee \otimes \det E \stackrel{\cong}\too E^{\vee\vee},
\end{equation}
a Lie sheaf becomes a tensor of the form
\begin{equation}
\label{tensor}
(F,\, \psi:F^{\otimes {r+1}}\too \SO_X),
\end{equation}
with $E=F\otimes \det F$. 

\begin{definition}
\label{lietensor}
A Lie tensor is a tensors of type $a=r+1$ which 
satisfies the following properties
\begin{enumerate}
\item $\psi$ factors through $F\otimes F\otimes \bigwedge^{r-1} F$,

\item the homomorphism $\wt \psi: F\otimes F \to F^{\vee\vee} \otimes 
\det F^\vee$ 
associated by (\ref{wedge}) 
is skew-symmetric.

\item the homomorphism $\wt \psi$ satisfies the Jacobi identity.
\end{enumerate}
\end{definition}

There is a canonical bijection between the set of isomorphism
classes of Lie sheaves $(E,\varphi:E\otimes E\to E^{\vee\vee})$
and Lie tensors $(F,\psi:F^{\otimes {r+1}}\too \SO_X)$
(with $E=F\otimes \det F$).

If the Lie algebra on the fiber $E(x)$ for all $x$ where
$E$ is locally free is always isomorphic
to a fixed semisimple Lie algebra $\fg$, then we say that it is 
a $\fg$-sheaf. Then, the Killing form gives an orthogonal
structure $\kappa:E\otimes E\to \SO_X$ to $E$.

\begin{definition}
A $\fg$-sheaf is (semi)stable if for all
orthogonal algebra filtrations, that is, filtrations with
$$
(1)\quad E_i^\perp = E_{-i-1}^{}
\quad \text{and} \quad
(2)\quad [E_i,E_i] \subset E_{i+j}^{\quad\vee\vee} 
$$
for all $i$, $j$, the following holds
$$
\sum (r P_{E_i} - r_i P_E ) (\leq) 0 \; .
$$
\end{definition}

It is shown in \cite{G-S2} that a $\fg$-sheaf is (semi)stable if and only 
if the associated tensor is $\delta$-(semi)stable, when 
$\delta$ has degree $n-1$.

We will sketch how the moduli space of tensors is constructed.
The idea is similar to the construction of the moduli space of torsion
free sheaves. First we construct a scheme which classifies
\textit{$\delta$-semistable based tensors}, 
that is, triples $(f,E,\varphi)$ where  
$f:V \to H^0(E(m))$ is an isomorphism, up to a constant, and
$(E,\varphi)$ is a $\delta$-semistable tensor. There is a natural
embedding of this scheme in a product $\PP(V_1)\times \PP(V_2)$, where
$V_1$ and $V_2$ are representations of $\slv$. An ample line bundle
with a linearization of the $\slv$ action is given by $\SO_X(b_1,b_2)$.
The choice of the integers $b_1$ and $b_2$ will depend on the polynomial
$\delta$, and the moduli space of $\delta$-semistable tensors will
be the GIT quotient.

%We already saw how a based sheaf $(f,E)$ gives a point in 
%$\PP(V_1)$, with $V_1=\bigwedge{}^{P(l)}(V^\vee\otimes W^\vee)$.
To find $V_2$, note that the isomorphism $f:V \to H^0(E(m))$
and $\varphi$ induces a linear map
$$
\Phi: V^{\otimes a} \;\too\; 
H^0(E(m)^{\otimes a}) \;\too\;
H^0(\SO_X(am)) \; =:\; B\; .
$$
Therefore, the semistable based tensor $(f,E,\varphi)$ gives a point
$(q,[\Phi])$
$$
\SH\times \PP(V_2)\; := \;
\Quot_X(V\otimes_\CC \SO_X(-m) ,P) \times
\PP\big((V^{\otimes a})^\vee \otimes B\big)
$$
The points obtained in this way have the property that the
homomorphism $\Phi$ composed with evaluation factors as
\[
\xymatrix{
{V^{\otimes a}\otimes \SO_X(-am)} \ar[d]^{\Phi\qquad} 
\ar@{->>}[r]^{q^{\otimes a}} &
{E^{\otimes a} \ar@(d,r)[ddl]^{\varphi}} \\ 
{H^0(\SO_X(am))\otimes \SO_X(-am)} 
\ar[d]^{\ev} & \\
{\SO_X}
}
\]
Let $Z'$ be the closed subscheme of $\SH\times \PP(V_2)$ where there
is a factorization as above, and let $Z\subset Z'$ be the closure
of the open subset $U\subset Z'$ of points 
$(q:V\otimes\SO_X(-m)\to E,[\Phi])$ such that the tensor is $\delta$-semistable.
Using Grothendieck's embedding $\SH\to \PP(V_1)$,
explained in section \ref{secsheaves}, we obtain a closed embedding
$$
Z \too \PP(V_1)\times \PP(V_2)
$$
We endow $Z$ with the polarization $\SO_Z(b_1,b_1)$, where
$$
\frac{b_2}{b_1}=\frac{P(l)\delta(m)-\delta(l)P(m)}{P(m)-a\delta(m)}
$$
In other words, we use the Segre embedding
$$
\PP(V_1)\times \PP(V_2) \too \PP(V_1^{\otimes b_1} \otimes
V_2^{\otimes b_2})
$$
and take the pullback of the ample line bundle $\SO_\PP(1)$.

It is proved in \cite{G-S1} that a point in $Z$ is GIT-(semi)stable if
and only if 
the induced linear map $f:V \to H^0(E(m))$ is an isomorphism
and it corresponds to a $\delta$-(semi)stable based tensor.
Therefore, the GIT quotient $Z\gitq \slv$ is the moduli space of
$\delta$-semistable tensors.

To show how this is used to obtain moduli spaces of related objects,
we will sketch the construction of the moduli space of orthogonal
sheaves. First we construct the projective scheme $Z$ as before,
for tensors of type $a=2$, i.e., of the form $(E,\varphi:E\otimes E
\to \SO_X)$. 
The condition of being symmetric is closed, so it
defines a closed subscheme $R\subset Z^{ss}$, and the GIT quotient
$R\gitq \slv$ is projective. On the other hand, 
the condition of being nondegenerate is open, so it defines
an open subscheme $R_1\subset R$. How can we prove
that, after we remove the points corresponding to degenerate bilinear
forms, the quotient is still projective?. The idea is to show
that, if $(E,\varphi)$ is degenerate, then it is $\delta$-unstable
(we remark that, to prove this, we need the degree of $\delta$
to be $n-1$).
Therefore, $R_1=R$, because all tensors 
corresponding to points in $R$ are semistable.

In other words, the moduli space of 
orthogonal sheaves $R_1\gitq \slv$ is projective because
the inclusion $R_1\inj R$ is proper (in fact, it is the identity).
Every time we impose a condition which is not
closed, we have to prove a properness result of this sort,
in order to show that the moduli space is projective.

The tensors defined in this section can easily be generalized
to tensors of type $(a,b,c)$, that is, pairs $(E,\varphi)$ consisting
of a torsion free sheaf and a homomorphism
\begin{equation}
\label{tensorabc}
\varphi: (E^{\otimes a})^{\otimes b}\too (\det E)^{\otimes c} \, .
\end{equation}
This more general notion will be needed in section \ref{secrho}.

\section{Principal bundles}

Recall that, in the \'etale topology, 
an open covering of a scheme $Y$ is a finite collection 
of morphisms $\{f_i:U_i\to Y\}_{i\in I}$ such that
each $f_i$ is \'etale, and $Y$ is the
union of the images of the $f_i$.

Note that an ``open \'etale subset'' of a scheme $Y$ 
is not really a subset of $Y$, but an \'etale morphism $U\to Y$.
If $f:X\to Y$ is a morphism, by a slight abuse of language
we will denote by $f^{-1}(U)$ the pull-back
\[
\xymatrix{
{f^{-1}(U)} \ar[r] \ar[d] & {X}\ar[d]^{f}\\
{U} \ar[r] & {Y}
}
\]
Let $G$ be an algebraic group. A principal $G$-bundle on $X$ is
a scheme $P$ with a right $G$-action and an invariant morphism
$P\to X$ with a $G$-torsor structure. 
A $G$-torsor structure is given by an atlas consisting on an
\'etale open covering $\{U_i\}$
of $X$ and $G$-equivariant isomorphisms $\psi_i:p^{-1}(U_i) \to
U_i\times G$, with $p=p_{U_i}\circ \psi_i$ (the 
$G$-action on $U_i\times G$ is given by 
multiplication on the right). Two atlases give the same $G$-torsor
structure if their union is an atlas.
An isomorphism of principal bundles is a $G$-equivariant 
isomorphism $\varphi:P \to
P'$. 

In short, a principal bundle is locally trivial in the \'etale
topology, and the fibers are $G$-torsors.
We remark that, if we were working in arbitrary characteristic,
an algebraic group could be non-reduced, and we should have
used the flat topology.

Given a principal $G$-bundle as above, we obtain an element of 
the \'etale cohomology set $\check H^1_\et(X,\underline{G})$,
and this gives a bijection between isomorphism classes of
principal $G$-bundles and elements of this set.
Indeed, since the isomorphisms $\psi_i$ of an atlas are required 
to be $G$-invariants, the composition $\psi^{}_j\circ\psi^{-1}_i$
is of the form $(x,g)\mapsto (x,\alpha_{ij}
(x)g)$, where $\alpha_{ij}:U_i\cap U_j\to G$ is
a morphism, which satisfies the cocycle condition and defines
a class in $\check H^1_\et(X,G)$.

Given a principal $G$-bundle $P\to X$ and a left 
action $\sigma$ of $G$ in
a scheme $F$, we denote
$$
P(\sigma,F) \;:=\; P\times_G F \;=\; (P\times F)/G,
$$
the associated fiber bundle.
Sometimes this notation is shortened to $P(F)$ or $P(\sigma)$.
In particular, for a representation $\rho$ of $G$ in
a vector space $V$, $P(V)$ is a vector bundle on $X$,
and if $\chi$ is a character of $G$, $P(\chi)$ is a line bundle.

If $\rho:G\to H$ is a group homomorphism, let $\sigma$ be
the action of $G$ on $H$ defined by left multiplication
$h\mapsto \rho(g)h$. Then, the associated fiber bundle
is a principal $H$-bundle, and it is denoted
$\rho_*P$. We say that this principal $H$-bundle is obtained
by \textit{extension of structure group}. 

Let $\rho:H \to G$ be a homomorphism of groups, and let $P$ be a
principal $G$-bundle on a scheme $Y$. 
A \textit{reduction of structure group} of $P$ 
to $H$ is a pair $(P^H,\zeta)$, where $P^H$ is 
a principal $H$-bundle on $Y$ and $\zeta$ is an isomorphism
between $\rho_* P^H$ and $P$. 
Two reductions $(P^{H},\zeta^{})$
and $(Q^{H},\theta^{})$ are isomorphic if there is an isomorphism
$\alpha$ giving a commutative diagram
\begin{equation}
\label{isopair}
\xymatrix{{P^H} \ar[d]^{\alpha}_{\cong} \\ 
{Q^{H}}  }
\qquad
\xymatrix{
  {\rho_{*} P^H} \ar[r]^{\zeta^{}}
\ar[d]^{\rho_{*}{\alpha}} & {P^{}} \ar@{=}[d] \\
 {\rho_{*} Q^{H}} \ar[r]^{\theta^{}} & {P^{}} 
}
\end{equation}
The names ``extension'' and ``reduction'' come from the case
in which $\rho$ is injective, but note that these notions
are still defined if the homomorphism is not injective.

If $\rho$ is injective, giving a reduction is equivalent to
giving a section $\sigma$ of the associated fibration $P(G/H)$, where
$G/H$ is the quotient of $G$ by the right action of $H$. 
Indeed, such a section gives a reduction $P^H$ by pull-back
\[
\xymatrix{
{P^H} \ar[r]^{i} \ar[d] & {P} \ar[d]\\
{X}\ar[r]^{\sigma} \ar[r] & {P(G/H)}
}
\]
and the isomorphism $\zeta$ is induced by $i$. Conversely,
given a reduction $(P^H,\zeta)$, the isomorphism $\zeta$ induces 
an embedding $i:P^H\to P$, and the quotient by $H$ of this
morphism gives a section $\sigma$ as above.

For example, if $G=\op{O}(r)$ and $H=\gl_r$, the quotient
$H/G$ is the set of non-degenerate bilinear symmetric forms
on the vector space $\CC^r$, hence a section of $P(H/G)$ is
just a non-degenerate bilinear symmetric morphism 
$E\otimes E\to \SO_X$, where $E$ is the vector bundle associated to
the principal $\gl_r$-bundle. 

To construct the moduli space, we will assume that $G$ is a connected
reductive algebraic group. Let $G'=[G,G]$ be the commutator subgroup,
and let $\fg=\fz\oplus \fgp$ be the Lie algebra of $G$, where 
$\fg'$ is the semisimple part and $\fz$ is the center.

Recall that, in the case of vector bundles, to obtain a projective
moduli space when $\dim X>1$, we had to consider also torsion free
sheaves. Analogously, principal $G$-bundles are not enough if we
want a projective moduli space, and this is why we also consider
principal $G$-sheaves, which we will now define.

\begin{definition}
\label{numero}
A principal $G$-sheaf $\SP$ over $X$ is a triple $\SP=(P,E,\psi)$
consisting of a torsion free sheaf $E$ on $X$, 
a principal $G$-bundle $P$ on the maximal open set $U_E$ where
$E$ is locally free, and an isomorphism of vector bundles
$$
\psi:P(\fgp)\stackrel{\isom}{\too} E|_{U_E} \, .
$$
\end{definition}

This definition can be understood from two points of view.
From the first point of view,
we have a torsion free sheaf $E$ on $X$, together with a reduction to $G$,
on the open set $U_E$, of the principal $\gl_r$-bundle corresponding
to the vector bundle $E|_{U_E}$.
Indeed, the pair $(P,\psi)$ is the same thing as a reduction
to $G$ of the principal $\gl_r$-bundle on $U_E$ 
associated to the vector bundle $E|_{U_E}$. 
It can be shown that, if we are given a reduction to a principal
$G$-bundle on a big open set $U'\subsetneq U_E$, this reduction
can uniquely be extended to $U_E$.

From the other point of view, we have a principal $G$-bundle on a big open
set $U$, hence a vector bundle $P(\fg')$, together with a given
extension of this 
vector bundle on $U$ to a torsion free sheaf on the whole of $X$.

The Lie algebra structure of $\fgp$ is semisimple, hence 
the Killing form is non-degenerate.
Correspondingly, the adjoint vector bundle
$P(\fgp)$ on $U$ has a Lie algebra structure
$P(\fgp)\otimes P(\fgp)\to P(\fgp)$ and an orthogonal
structure, $\kappa:P(\fgp)\otimes P(\fgp)\to \SO_U$.
These uniquely extend to give orthogonal and
$\fgp$-sheaf structure to $E$:
$$
\kappa:E\otimes E \too \SO_X \qquad
[\,,] : E\otimes E \too E^{\vee\vee}
$$
where we have to take $E^{\vee\vee}$ in the target
because an extension $E\otimes E\to E$ does not
always exist. 
The orthogonal structure assigns an orthogonal 
$F^\perp=\ker(E\inj E^\vee \to F^\vee)$
to each subsheaf $F\subset E$.

\begin{definition}
\label{stab1}
A principal $G$-sheaf $\SP=(P,E,\psi)$ is said to be 
(semi)stable if for 
all orthogonal algebra filtrations 
$E_\bullet\subset E$, that is, filtrations with 
$$
(1)\quad E_i^\perp = E_{-i-1}^{}
\quad \text{and} \quad
(2)\quad [E_i,E_i] \subset E_{i+j}^{\quad\vee\vee} 
$$
for all $i$, $j$, the following holds
$$
\sum (r P_{E_i} - r_i P_E ) (\leq) 0
$$
\end{definition}

Replacing the Hilbert polynomials $P_E$ and $P_{E_i}$ by degrees, 
we obtain the notion of \textit{slope (semi)-stability}.

Clearly 
$$
\text{slope-stable $\Longrightarrow$ stable
$\Longrightarrow$ semistable $\Longrightarrow$ 
slope-semistable}
$$

Since $G/G'\isom
\CC^{*q}$, given a principal $G$-sheaf, the principal
bundle $P(G/G')$ obtained by extension of structure group
provides $q$ line bundles on $U$, and since 
$\codim X\setminus U \geq 2$, 
these line bundles extend uniquely to line bundles on $X$.
Let $d_1,\ldots,d_q\in H^2(X,\CC)$ be their Chern classes. 
The rank $r$ of $E$ is clearly the dimension of $\fgp$.
Let $c_i$ be the Chern classes of $E$.

\begin{definition}[Numerical invariants]
\label{ginvariants}
We call the 
data $\tau=(d_1,\ldots,d_q,c_i)$ 
the numerical invariants of the principal
$G$-sheaf $(P,E,\psi)$. 
\end{definition}

\begin{definition}[Family of semistable principal $G$-sheaves]
A family of (semi)stable principal $G$-sheaves 
parameterized by a scheme $S$
is a triple 
$(P_S,E_S,\psi_S)$, 
with $E_S$ a family of torsion free sheaves,
$P_S$ a principal $G$-bundle on 
the open set $U_{E_S}$ where $E_S$ is locally free, 
and $\psi:P_S(\fgp)\to E_S|_{U_{E_S}}$ an isomorphism of
vector bundles, such that for
all closed points $s\in S$ the corresponding principal $G$-sheaf
is (semi)stable with numerical invariants $\tau$.
\end{definition}

An isomorphism between two such families $(P_S,E_S,\psi_S)$ and
$(P'_S,E'_S,\psi'_S)$ is a pair 
$$
(\beta:P^{}_S\stackrel{\isom}\too P'_S ,\gamma:E^{}_S \stackrel{\isom}\too E'_S)
$$
such that the following diagram is commutative
$$
\xymatrix{
{P_S(\fgp)} \ar[r]^{\psi} \ar[d]_{\beta(\fgp)} & 
{E_S|_{U_{E_S}}} \ar[d]^{\gamma|_{U_{E_S}}} \\
{P'_S(\fgp)} \ar[r]^{\psi'} & {E'_S|_{U_{E_S}}} 
}
$$
where $\beta(\fgp)$ is the isomorphism of vector bundles induced by
$\beta$.
Given an $S$-family $\SP_S=(P_S,E_S,\psi_S)$ and a morphism
$f:S'\to S$, the pullback is defined as
$\wt{f}^* \SP_S=(\wt{f}^*P_S ,\overline{f}^* E_S,\wt{f}^*\psi_S)$,
where $\overline{f}=\id_X\times f:X\times S  \to 
X\times S'$ and
$\wt{f}=i^*(\overline{f}):U_{\overline{f}^*E_S}\to U_{E_S}$,
denoting $i:U_{E_S}\to X\times S$ the inclusion of the open
set where $E_S$ is locally free.

We can then define the functor of families of semistable 
principal $G$-sheaves 
$$
F^\tau_G: (\Sch/\CC) \too (\Sets) 
$$
sending a scheme $S$, locally of finite type, 
to the set of isomorphism classes of families
of semistable principal 
$G$-sheaves with numerical invariants $\tau$. As usual, it is
defined on morphisms as pullback.

\begin{theorem}
There is a projective
moduli space of semistable $G$-sheaves on $X$
with fixed numerical invariants.
\end{theorem}

This theorem is a generalization of the theorem of Ramanathan,
asserting the existence of a moduli space of semistable principal 
bundles on a curve. 

Note that in the definition of principal $G$-sheaf we have used
the adjoint representation on the semisimple part $\fg'$ of the Lie
algebra of $G$, to obtain a vector bundle $P(\fg')$ on a big open set
of $X$, which we extend to the whole of $X$ by torsion free sheaf. 
If we use a different representation $\rho:G\to \gl_r$, we have the
notion of principal $\rho$-sheaf:

\begin{definition}
A principal $\rho$-sheaf $\SP$ over $X$ is a triple $\SP=(P,E,\psi)$
consisting of a torsion free sheaf $E$ on $X$, 
a principal $G$-bundle $P$ on the maximal open set $U_E$ where
$E$ is locally free, and an isomorphism of vector bundles
$$
\psi:P(\rho)\stackrel{\isom}{\too} E|_{U_E} \, .
$$
\end{definition}

Now we will give some examples of principal $\rho$-sheaves which have
already appeared:

\begin{itemize}

\item If $G=\gl_r$ and $\rho$ is the canonical representation,
then a principal $\rho$-sheaf is a torsion free sheaf.

\item If $G=\orth$ and $\rho$ is the canonical representation,
then a principal $\rho$-sheaf is an orthogonal sheaf.

\item If $G=\sor$ and $\rho$ is the canonical representation,
then a principal $\rho$-sheaf is a special orthogonal sheaf 
(cf. \cite{G-S1}), that is, a triple $(E,\varphi,\psi)$ where
$\varphi:E\otimes E\to \SO_X$ symmetric and nondegenerate, and
$\psi:\det E\to \SO_X$ is an isomorphism such that $\det
\varphi=\psi^{\otimes 2}$.

\item If $G=\sympl$ and $\rho$ is the canonical representation,
then a principal $\rho$-sheaf is a symplectic sheaf.

\item If $G$ is semisimple and $\rho$ is injective, then
giving a principal $\rho$-sheaf is equivalent to giving 
a honest singular principal bundle \cite{Sch1,Sch2} with respect
to the dual representation $\rho^\vee$ (see section \ref{secrho}).

\end{itemize}

In all these cases (and also for principal $G$-sheaves, i.e.,
when $\rho:G\to \gl(\fg')$ is the adjoint
representation), the stability condition is equivalent to
the following:

\begin{definition}[Stability for principal $\rho$-sheaves]
A principal $\rho$-sheaf $\SP=(P,E,\psi)$ is said to be 
(semi)stable if for all reductions on any big open set $U\subset U_E$ 
of $P$ to a parabolic subgroup $Q\subsetneq G$, and all dominant
characters of $Q$, which are trivial on the center of $Q$, the
induced filtration of saturated torsion free sheaves
$$
\ldots \subset E_{i-1}\subset E_{i}
\subset E_{i+1}\subset \ldots
$$
satisfies the following
$$
\sum (r P_{E_i} - r_i P_E ) (\leq) 0
$$
\end{definition}

\section{Construction of the moduli space of principal sheaves}

In this section we will give a sketch of the construction of the
moduli space in \cite{G-S2}. The strategy is close to that of
Ramanathan.

Let $r=\dim \fgp$, and consider the adjoint representation 
$\rho:G \to \gl_r$  of $G$ in $\fgp$. 
The idea of Ramanathan is to start by
constructing a scheme $R_0$ which classifies based vector bundles
of rank $r$, and then to construct another scheme $Q\to R_0$ 
such that
the fiber over each based vector bundle $(f,E)$ parameterizes all
reductions to $G$ of the principal $\gl_r$-bundle $E$. In other
words, $Q$ classifies tuples $(f,P,E,\psi)$, where $f$ is an
isomorphism of a fixed vector space $V$ with $H^0(E(m))$,
$P$ is a principal $G$-bundle and $\psi$ is an isomorphism between
the vector bundle $P(\rho,\fgp)$ and $E$.

The problem is that $\rho$ is not injective in general, so it is
not easy to construct a reduction of structure group from
$\gl_r$ to $G$ in one step. 
Therefore, Ramanathan
factors the representation $\rho$ into several 
group homomorphisms,
and then constructs reductions step by step. 

Recall that $G'=[G,G]$ is
the commutator subgroup. Let $Z$ (respectively, $Z'$) 
be the center of $G$ (respectively, $G'$). Note that $Z'=G'\cap Z$.
The adjoint representation factors as follows
\[
\xymatrix{
{G} \ar@{->>}^{\rho_3}[r]&
{G/Z'}  
\ar@{->>}^{\rho'_2}[r]&
{G/Z} 
\ar@{^{(}->}^{\rho_2}[r]&
{\Aut(\fgp)} 
\ar@{^{(}->}^{\rho_1}[r]&
{\gl_r}
}
\]
and the schemes parameterizing these reductions are
$$
R_3 \stackrel{f_3}\too R'_2 \stackrel{f'_2}\too R_2 
\stackrel{f^{}_2}\too R_1\too R_0
$$
In the case $\dim X=1$ this works well because a principal $G$-bundle
is semistable if and only if the associated vector bundle is semistable.
This is no longer true if $X$ is not a curve,
and this is why, for arbitrary dimension, 
we do not construct the scheme $R_0$, but instead start directly
with a scheme $R_1$, classifying 
semistable based principal $\Aut(\fgp)$-sheaves.

Here $\Aut(\fgp)$ denotes the subgroup of $\gl_r$ 
of linear automorphisms which respect the Lie algebra structure.
Therefore, a based principal $\Aut(\fgp)$-sheaf  is the
same thing as a based $\fgp$-sheaf. 

Using the isomorphism (\ref{wedge}), we can describe a $\fgp$-sheaf
as a Lie tensor (definition \ref{lietensor}).
such that the Lie algebra structure induced on the fibers of $E$, over
points $x\in X$ where $E$ is locally free, is isomorphic to $\fgp$.

Choose a polynomial $\delta$ of degree $\dim X-1$, with positive leading
coefficient. We fix the first Chern class to be zero. This is because
we are interested in $\fgp$-sheaves, and since $\fgp$ is semisimple,
its Killing form is nondegenerate, hence induces an orthogonal
structure on the sheaf, and this forces the first Chern class to
be zero.

We start with the scheme $Z$, defined in section
\ref{sectensors}, classifying based tensors of
type $a=r+1$. This scheme has an
open subset $Z^{ss}$ corresponding to $\delta$-semistable tensors.
Conditions (1) to (3) in the definition of Lie tensor are closed, 
hence they define a closed subscheme
$R\subset Z^{ss}$. Using the isomorphism (\ref{wedge}), we
see that the scheme $R$ parameterizes Lie sheaves. Recall
that a Lie sheaf structure 
induces a Killing form $\kappa:E\otimes E\to \SO_X$.

\begin{lemma}
There is a subscheme $R_1\subset R$ corresponding to those
Lie tensors which are $\fgp$-tensors.
\end{lemma}

The family of Lie sheaves parameterized by $R$ gives a family
of Killing forms $E_R\otimes E_R\to \SO_{X\times R}$, and
hence a homomorphism $f:\det E_R\to \det E_R^\vee$. We have
fixed the determinant of the tensors to be trivial, hence
$\det E_R$ is the pullback of a line bundle on $R$,
and therefore the homomorphism $f$ is nonzero on an open 
set of the form $X\times W$, where $W$ is an open set of
$R$. The open set $W$ is in fact the whole of $R$. This
is because if $z$ is a point in the complement, it corresponds
to a Lie sheaf whose Killing form is non-degenerate, and hence
has a nontrivial kernel. Using this, it is possible to 
construct a filtration which shows that this
Lie sheaf is $\delta$-unstable when $\deg \delta=\dim X-1$, 
but this contradicts
the fact that $R\subset Y^{ss}$.

The Killing form of a Lie algebra is semisimple if and only if
it is non-degenerate. Therefore, for all points $(x,t)$ in the
open subset $\SU_{E_R}\subset X\times R$ where $E_R$ is locally free,
the Lie algebra is semisimple. 

Semisimple Lie algebras are rigid, that is, if there is a family
of Lie algebras, the subset of the parameter space corresponding
to Lie algebras isomorphic to a given semisimple Lie algebra is
open. Therefore, since $U_{E}$ is connected for all torsion free
sheaves $E$, all points $(x,t)\in \SU_{E_R}\subset X\times R$ 
where $t$ is in a fixed connected component of $R$, give isomorphic
Lie algebras. Let $R_1$ be the union of those components whose Lie
algebra is isomorphic to $\fgp$. The inclusion 
$$
i:R_1 \inj R
$$
is proper, and hence, since the GIT quotient $R\gitq \slv$ is
proper, also the GIT quotient $R_1\gitq \slv$ is proper.
Note that, to prove properness of $i$, two facts about
semisimple Lie algebras were used: rigidity, and 
nondegeneracy of their Killing forms.

For simplicity of the exposition, to explain the successive
reductions, first we will assume that for 
all $\fgp$-sheaves $(E,\varphi)$, the torsion free sheaf $E$ is 
locally free. In other words, $U_E=X$ (this holds, for instance,
if $\dim X=1$). At the end we will mention what has to be modified
in order to consider the general case.

The group $G/Z$ is the connected component of identity of
$\Aut(\fgp)$. 
Therefore, giving a reduction of structure 
group of a principal $\Aut(\fgp)$-bundle $P$
by $\rho_2$ is the same thing as giving a section
of the finite \'etale morphism $P(F)\to X$, 
where $F$ is the finite group $\Aut(\fgp)/(G/Z)$.
This implies that $R_2\to R_1$ is a finite \'etale morphism,
whose image is a union of connected components of $R_1$.

There is an isomorphism of groups 
$G/Z'\cong G/G' \times G/Z$, and $\rho'_2$ is just the projection
to the second factor. Therefore, a reduction to $G/Z'$ of a principal
$G/Z$-bundle $P^{G/Z}$ is just a pair $(P^{G/G'},P^{G/Z})$, where 
$P^{G/Z}$ is the original $G/Z$-bundle and $P^{G/G'}$ is a
$G/G'$-bundle. But 
$$G/G'\cong \CC^*\overbrace{\times\cdots\times}^q \CC^* \;,$$
hence this is just a collection of $q$ line bundles,
whose Chern classes are given by the numerical invariants
which have been fixed. 
This implies that there is an isomorphism
$$
R'_2 \cong J\overbrace{\times\cdots\times}^q J \times R_2 \, ,
$$
where $J$ is the Jacobian of $X$.

Finally, we have to consider reductions of a principal
$G/Z'$-bundle to $G$, where $Z'$ is a finite subgroup of the center of $G$.
There is an exact sequence of pointed sets (the distinguished point 
being the trivial bundle)
$$
\check H^1_\et(X,\underline{Z'}) \too
\check H^1_\et(X,\underline{G}) \too
\check H^1_\et(X,\underline{G/Z'}) \stackrel{\delta}\too
\check H^2_\et(X,\underline{Z'}) \, .
$$
Note that $Z'$ is abelian, therefore $H^i_\et(X,\underline{Z'})$ is
an abelian group, and it is isomorphic to the singular cohomology
group $H^i(X;Z')$, hence finite.
A principal $G/Z'$-bundle admits a reduction to $G$ if and only
if the image by $\delta$ of the corresponding point is $0$.
This is an open and closed condition, therefore there is a subscheme 
$\hat R'_2$ of $R'_2$, consisting of a union of connected components,
corresponding to those principal $G/Z$-bundles admitting a reduction
to $G$.

Let $(P^G,\zeta)$ be a reduction to $G$ of a principal $G/Z'$-bundle.
It can be shown that the set of isomorphism classes of all 
reductions to $G$ is in bijection with the cohomology set
$
\check H^1_\et(X,\underline{Z'}),
$
with the unit element of this set corresponding to the chosen
reduction $(P^G,\zeta)$. This cohomology set is an abelian group, 
because $Z'$ is abelian. Therefore, the set of reductions of
a principal $G/Z'$-bundle to $G$ is an $\check
H^1_\et(X,\underline{Z'})$-torsor, and this implies that
$R_3\to \hat R_2'$ is a principal $\check
H^1_\et(X,\underline{Z'})$-bundle. 
Using that this cohomology set is 
a finite set (in fact isomorphic to the singular cohomology
group $H^1(X;Z')$), and that $\hat R'_2$ is a union of connected
components of $R'_2$, it follows that
$R_3\to R'_2$ is finite \'etale.

Ramanathan \cite[Lemma 5.1]{Ra} proves that, 
if $H$ is a reductive algebraic group, 
$f:Y\to S$ is an $H$-equivariant affine morphism, and $p:S\to
\overline{S}$ is a good quotient, then $Y$ has a good quotient
$q:\overline{Y}\to Y$ and the induced morphism $\overline{f}$
is affine. Moreover, if $f$ is finite, $\overline{f}$ is also
finite. When $f$ is finite and $p$ is a geometric quotient, 
also $q$ is a geometric quotient.

The group $\slv$ acts on all the schemes $R_i$, and the morphisms
$f_2$ and $f_3$
are equivariant and finite.
Therefore, we can apply Ramanathan's lemma to those morphisms.

The morphism $f'_2:J^{\times q}\times R_2\to R_2$ 
is just projection to a factor, and the group acts trivially on 
the fiber, therefore if $p_2:R_2\to \FM_2$ is a good quotient of $R_2$,
$J^{\times q}\times \FM_2$ will give a good quotient of $R'_2$.
Furthermore, if $p_2$ becomes a geometric quotient when restricting
to an open set, the same will be true after taking the product with
$J^{\times q}$.

Using GIT, we know that $R_1$ has a good quotient $\FM_1$, which is 
a geometric quotient when restricting to the open set of stable points.
Therefore, the same holds for all these schemes, and the good quotient
of $R_3$ is the moduli space of principal $G$-bundles.

%A key point of Ramanathan construction is that a principal $G$-bundle
%is semistable if and only if the associated vector bundle $E=P(\fgp)$
%is semistable as a vector bundle 
%(incidentally, it can happen that a stable principal
%$G$-bundle gives a vector bundle which is strictly semistable).
%
%If $\dim X>1$ this is not true. What is true in higher dimensions
%is that a principal $G$-sheaf is semistable if and only if the
%associated Lie algebra sheaf is semistable (as a Lie algebra sheaf).
%And this is why, in higher dimensions,
%the starting point is not $R$, but $R_1$, the parameter space of 
%based Lie algebra sheaves.

The successive reductions in higher dimension are very similar to
the reductions in the case $X$ is a curve, except for the technical
difficulty that the principal bundles in general are not defined in
the whole of $X$, but only in a big open set. To overcome this
difficulty, we need ``purity'' results for open sets $U\subset X$ when
$U$ is big. We will discuss them one by one.

First we consider reductions
of a principal $\Aut(\fgp)$-bundle $P$
to $G/Z$. These are parameterized by sections of the associated
fibration $P(F)$, where $F=\Aut(\fgp)/(G/Z)$ is a finite group.
If $P$ is a principal bundle on a big open set $U$, $P(F)$ is a
Galois cover of $U$, given by a representation of the algebraic
fundamental group of $\pi(U)$ in $F$. Since $U$ is a big open set, 
$\pi(X)=\pi(U)$ (purity of fundamental group), and hence the
Galois cover $P(F)$ of $U$ extends uniquely to a Galois cover
of $X$. This implies that, even if $\dim X>1$, the morphism
$R_2\to R_1$ is still finite \'etale, as in the curve case.

Giving a reduction of a principal $G/Z$-bundle 
on $U$ to a principal $G/Z'$-bundle is equivalent to giving $q$ line
bundles on $U$. Since $U$ is a big open set, the Jacobians of
$U$ and $X$ are isomorphic (purity of Jacobian), and hence
we still have $R'_2=J(X)^{\times q} \times  R_2$.

Finally, we have to consider reductions of principal $G/Z'$-bundles
to $G$. Using the fact that $U$ is a big open set, there
are isomorphisms $\check H^i_\et(X,\underline{Z'})
\cong \check H^i_\et(U,\underline{Z'})$ for $i=1,2$.
Therefore, the arguments used for the case $U=X$ still hold
in general, and it follows that $R_3\to R'_2$ is \'etale finite.

\section{Construction of the moduli space of principal $\rho$-sheaves}
\label{secrho}

In \cite{Sch1,Sch2}, A. Schmitt fixes a semisimple group $G$ and
a faithful representation $\rho$, defines semisimple honest singular
principal bundle with respect to this data (see definition below), 
and constructs the corresponding
projective moduli space. 
Giving such an object is equivalent to giving a principal
$\rho^\vee$-bundle, where $\rho^\vee$ is the dual representation in 
$V^\vee$. In this section we will give a sketch of Schmitt's construction.

Let $G$ be a semisimple group, and 
$\rho:G\to \glv$ a faithful representation.
A honest singular principal $G$-bundle is a pair $(\SA,\tau)$,
where $\SA$ is a torsion free sheaf on $X$ and
$$
\tau: \Sym^*(\SA\otimes V)^G \too \SO_X
$$
is a homomorphism of $\SO_X$-algebras
such that, if  
$
\sigma: X \to Hom(V\otimes \SO_U, \SA|^\vee_U)\gitq G
$
is the induced morphism, then 
$$
\sigma(U) \subset  Isom(V\otimes \SO_U, \SA|^\vee_U)/G \subset
Hom(V\otimes \SO_U, \SA|^\vee_U)\gitq G\, .
$$
It can be shown that the points in the affine $U$-scheme
$Isom(V\otimes_\CC \SO_U,\SA|^\vee_U)$
are in the open set of GIT-polystable points of 
$Hom(V\otimes_\CC \SO_U,\SA|^\vee_U)$, under the natural action of $G$,
therefore the previous inclusion makes sense.

Note that the homomorphism $\tau$ is uniquely defined by its 
restriction to $U\subset X$, therefore, 
giving a honest singular principal $G$-bundle is equivalent
to giving a principal $\rho^\vee$-sheaf $(P,E,\psi)$,
where 
$\rho^\vee:G\to \gl(V^\vee)$ is the dual representation,
$P$ is a principal $\gl_n$-bundle, $E=\SA$, and $\psi$
is induced by $\sigma|_U$.

In other words, in a principal $\rho$-sheaf, we extend to the whole
of $X$, as a
torsion free sheaf $E$, the vector bundle associated to $\rho$,
whereas, is a honest singular principal $G$-bundle associated
to $\rho$, we extend the dual of the vector bundle associated
to $\rho$.

The idea of Schmitt's construction is to transform $\tau$ into
a tensor.
Note that $\tau$ is an infinite collection
of $\SO_X$-module homomorphisms 
\begin{equation}
\label{taui}
\tau_i: \Sym^i(\SA\otimes V)^G \too \SO_X,
\end{equation}
but, since $\Sym^*(\SA\otimes V)^G$ is finitely generated as
a $\SO_X$-algebra, there is an integer $s$ such that
\begin{enumerate}
\item the sheaf
$$
\bigoplus_{i=1}^s \Sym^i(\SA\otimes V)^G
$$ 
contains a set of generators of the algebra, and
\item the subalgebra 
$$
\Sym^{(s!)}(\SA\otimes V)^G \;:=\; \bigoplus_{m=0}^\infty
\Sym^{s!m}(\SA\otimes V)^G
$$
is generated by elements in $\Sym^{s!}(\SA\otimes V)^G\,$.
\end{enumerate}
Using the homomorphisms $\tau_s$, we construct a homomorphism
of $\SO_X$-modules
\begin{equation}
\label{surj1}
\bigoplus_{\sum id_i=s!}
\Big(\bigotimes_{i=1}^s \Sym^{d_i}\big(\Sym^{i}(\SA\otimes V)^G\big)\Big)
\surj
{\Sym^{s!}(\SA\otimes V)^G}
\stackrel{\tau_s}\too
{\SO_X}
\end{equation}
Note that the vector space 
$$
\bigoplus_{\sum id_i=s!}
\Big(\bigotimes_{i=1}^s \Sym^{d_i}\big( \Sym^{i}(\CC^r\otimes V)^G\big)\Big)
$$
has a canonical representation of $\gl_n$, homogeneous of degree
$s!$, and hence it is a quotient of the representation
$$
(\CC^{\otimes a})^{\oplus b} \otimes (\bigwedge^r \CC^r)^{-\otimes c}
$$
for appropriate values of $a$, $b$ and $c$.
Therefore, there is a surjection
\begin{equation}
\label{surj2}
{(\SA^{\otimes a})^{\oplus b} \otimes (\det \SA)^{-\otimes c}}
\surj
{\bigoplus_{\sum id_i=s!}
\Big(\bigotimes_{i=1}^s \Sym^{d_i}\big(\Sym^{i}(\SA\otimes V)^G\big)\Big)}
\end{equation}
and composing (\ref{surj1}) with (\ref{surj2}) 
we obtain a tensor of type $(a,b,c)$, as in (\ref{tensorabc}).

\end{document}